\documentclass[12pt]{article}
\usepackage{amsmath}
\usepackage{latexsym}
\usepackage{amsfonts}
\usepackage[normalem]{ulem}
\usepackage{array}
\usepackage{amssymb}
\usepackage{graphicx}
\usepackage[backend=biber,
style=numeric,
sorting=none,
isbn=false,
doi=false,
url=false,
]{biblatex}

\usepackage{subfig}
\usepackage{wrapfig}
\usepackage{wasysym}
\usepackage{enumitem}
\usepackage{adjustbox}
\usepackage{ragged2e}
\usepackage[svgnames,table]{xcolor}
\usepackage{tikz}
\usepackage{longtable}
\usepackage{changepage}
\usepackage{setspace}
\usepackage{hhline}
\usepackage{multicol}
\usepackage{tabto}
\usepackage{float}
\usepackage{multirow}
\usepackage{makecell}
\usepackage{fancyhdr}
\usepackage[toc,page]{appendix}

\usetikzlibrary{shapes.symbols,shapes.geometric,shadows,arrows.meta}
\tikzset{>={Latex[width=1.5mm,length=2mm]}}
\usepackage{flowchart}\usepackage[paperheight=11.69in,paperwidth=8.27in,left=1.0in,right=1.0in,top=1.0in,bottom=1.0in,headheight=1in]{geometry}
\usepackage[utf8]{inputenc}
\usepackage[T1]{fontenc}
\TabPositions{0.5in,1.0in,1.5in,2.0in,2.5in,3.0in,3.5in,4.0in,4.5in,5.0in,5.5in,6.0in,}

\setlistdepth{9}
\renewlist{enumerate}{enumerate}{9}
		\setlist[enumerate,1]{label=\arabic*)}
		\setlist[enumerate,2]{label=\alph*)}
		\setlist[enumerate,3]{label=(\roman*)}
		\setlist[enumerate,4]{label=(\arabic*)}
		\setlist[enumerate,5]{label=(\Alph*)}
		\setlist[enumerate,6]{label=(\Roman*)}
		\setlist[enumerate,7]{label=\arabic*}
		\setlist[enumerate,8]{label=\alph*}
		\setlist[enumerate,9]{label=\roman*}

\renewlist{itemize}{itemize}{9}
		\setlist[itemize]{label=$\cdot$}
		\setlist[itemize,1]{label=\textbullet}
		\setlist[itemize,2]{label=$\circ$}
		\setlist[itemize,3]{label=$\ast$}
		\setlist[itemize,4]{label=$\dagger$}
		\setlist[itemize,5]{label=$\triangleright$}
		\setlist[itemize,6]{label=$\bigstar$}
		\setlist[itemize,7]{label=$\blacklozenge$}
		\setlist[itemize,8]{label=$\prime$}

\begin{document}
\setstretch{2.0}
\begin{justify}
{\fontsize{13pt}{15.6pt}\selectfont \textbf{A dual identity based symbolic understanding of the Gödel’s incompleteness theorems, P-NP problem, Zeno’s paradox and Continuum Hypothesis }\par}
\end{justify}\par

\begin{FlushRight}
\tab \tab \tab \tab \tab \tab \tab \tab {\fontsize{13pt}{15.6pt}\selectfont \textbf{\textit{Arun Uday}}\par}
\end{FlushRight}\par

\begin{justify}
\textit{Abstract: A semantic analysis of formal systems is undertaken, wherein the duality of their symbolic definition based on the $``$State of Doing$"$  and $``$State of Being$"$  is brought out. We demonstrate that when these states are defined in a way that opposes each other, it leads to contradictions. This results in the incompleteness of formal systems as captured in the Gödel’s theorems. We then proceed to resolve the P-NP problem, which we show to be a manifestation of Gödel’s theorem itself. We then discuss the Zeno’s paradox and relate it to the same aforementioned duality, but as pertaining to discrete and continuous spaces. We prove an important theorem regarding representations of irrational numbers in continuous space. We extend the result to touch upon the Continuum Hypothesis and present a new symbolic conceptualization of space, which can address both discrete and continuous requirements. We term this new mathematical framework as $``$hybrid space$"$ .}
\end{justify}\par

\begin{justify}
Keywords: Gödel, P-NP, incompleteness, complexity, Zeno, Continuum Hypothesis 
\end{justify}\par

\begin{enumerate}
	\item \textbf{Gödel’s incompleteness theorems }\par

Ever since Gödel came up with his famous theorems (Gödel [1931]) on incompleteness of formal systems, there have been numerous interpretations of them. Our goal here is not to undertake a full analysis of them, which has already been done by others (Detlefsen [1979], Smith [2007], Berto [2011], Kennedy [2014]). Rather, we attempt to build a semantic understanding of the theorems specifically based on the duality of identity states. This will enable us to establish the causal basis of the theorems. \par

We begin by stating that any formal system is \textit{defined} by two semantic identity states – the $``$State of Doing$"$  and the $``$State of Being$"$ . To cite a simple analog from common parlance, the \textit{identity} of a $``$driver$"$  is \textit{defined }as a person who drives. The \textit{identity} of a $``$teacher$"$  is \textit{defined} as a person who teaches, etc. In other words, their identity states are defined by their occupational states. Paradoxes occur when these two states, i.e. who they are and what they do contradict each other especially when their actions are self-referenced. This is best illustrated in the famous Russel’s paradox (Coffa [1979]), concerning the barber who shaves only those who do not shave themselves. In this case, the barber’s two states are:\par

\begin{itemize}
	\item \textit{State of Doing:} Shaving someone who does not shave himself\par

	\item \textit{State of Being:} Shaves someone who does not shave himself
\end{itemize}\par

When the $``$someone$"$  here is referenced to an external person who does not shave himself, there is no problem. But when the $``$someone$"$  is self-referenced to the barber himself, it causes a contradiction owing to the conflicting requirements of the States of Doing and Being. \par

Now, let us extend this understanding to formal systems. We know that a formal system is defined by its axioms and rules. These in turn define the set of theorems, which can be derived from them in the form of proofs of the theorems. Hence, in a similar manner as the examples above, for a formal system, its two symbolic states are:\par

\begin{itemize}
	\item \textit{State of Doing:} Proving a proposition (correctly)\par

	\item \textit{State of Being:} Proves a proposition (as being true)
\end{itemize}\par

As with the barber example, when the proposition itself is defined as something that is false and then self-referenced, it causes a contradiction of the kind \textit{P}, \par

\textit{P}:= $``$PROVE THIS PROPOSITION TO BE FALSE$"$ . \par

The contradiction arises because if the formal system does a correct job of the above task, its State of Doing will be – $``$Proving \textit{P} to be false$"$ . But if it successfully does that, it implies that $``$\textit{P} is false$"$  is false, or in other words \textit{P} is true and hence its State of Being will be - $``$Proves \textit{P} to be true$"$ . Clearly now, there is a contradiction between the system’s State of Doing and its State of Being. This is what we will term as mutually opposing dual symbolic states. Gödel basically formalized this in the form of his famous proof of incompleteness. Through the Gödel numbering system, he was able to exploit this duality of symbolic states to construct propositions of the kind \textit{P} above to demonstrate the incompleteness of any formal system. Since such propositions cannot be proven to be either true or false, the conclusion is that the formal system is incomplete. To be clear, not all dual symbolic identities result in paradoxes, but only mutually contradictory ones. The intuition underlying this duality will form a key aspect of our analysis, which we will carry through the rest of the paper. We will also demonstrate why such contradictions are impossible to avoid through any meta-level exclusions.\par

	\item \textbf{P-NP problem}\par

Ever since Cook formally introduced the P-NP problem (Cook [1971]), it has gone on to become one of the most fundamental unsolved problems in theoretical computer science. However, the history of the problem can be traced to a much earlier time, indeed to Gödel himself. It was discovered that in 1956, Gödel in a letter written to von Neumann enquired him about his views on the complexity of an NP complete problem (Hartmanis [1989)). As is well documented, the fundamental P-NP question is - if every problem whose solution can be verified in polynomial time can also be solved in polynomial time. We will now prove that it cannot be done so.\par

For this, let us construct a simple example called the $``$Alice Bob game$"$ . The game is seemingly an easy one. Alice thinks of a sequence of indefinite length of 0’s and 1’s, which Bob needs to start guessing. The moment he gets a character wrong, he loses\footnote{ The assumption here is that Alice won’t lie, i.e.- she won’t say that Bob is wrong even though he has correctly guessed the character that Alice had in her mind }. Bob agrees to play the game and they get started. Bob tries a few times and soon gives up as he is not able to go beyond guessing a few characters. But, he has a better idea! He writes a computer program, which will exactly predict the sequence of numbers. He is convinced of the brilliant logic that he has employed and now approaches Alice with great aplomb. Alice is however unperturbed. In fact, she goes on to prove that Bob’s computer program can never predict the exact sequence no matter what logic he has used. Her proof relies on the familiar Cantor diagonalization method. \par

Alice knows which programming language Bob uses for his programs (and it could be anything as it doesn’t affect the validity of the proof). She conceives of a list,  \( L \)  of all possible computer programs, which can be written using the language. This is just a subset of all the permutations of strings that can be formed using the vocabulary of the language. Since these can be ordered in an increasing order of their lengths, they will be countably infinite in number. A subset of these programs will be those that are able to generate sequences of 0’s and 1’s. Let us say, the ordered list  \( L^{'} \)  of these programs is P\textsubscript{1}, P\textsubscript{2}, P\textsubscript{3}, P\textsubscript{4} etc. \par

Clearly,  \( L^{'} \subset L \) \par

Now, she envisages the possible outputs that these programs can deliver. These are captured in the following table wherein C\textsubscript{1}, C\textsubscript{2}, C\textsubscript{3}, C\textsubscript{4} $ \ldots $  are the characters of the sequence generated by each of the programs P\textsubscript{1}, P\textsubscript{2}, P\textsubscript{3}, P\textsubscript{4 $ \ldots $ } \par

%%%%%%%%%%%%%%%%%%%% Table No: 1 starts here %%%%%%%%%%%%%%%%%%%%

\begin{table}[H]
 			\centering
\begin{tabular}{p{0.89in}p{0.29in}p{0.29in}p{0.28in}p{0.29in}p{2.71in}}
\hline
%row no:1
\multicolumn{1}{|p{0.89in}}{\textbf{Program}} & 
\multicolumn{1}{|p{0.29in}}{\textbf{C\textsubscript{1}}} & 
\multicolumn{1}{|p{0.29in}}{\textbf{C\textsubscript{2}}} & 
\multicolumn{1}{|p{0.28in}}{\textbf{C\textsubscript{3}}} & 
\multicolumn{1}{|p{0.29in}}{\textbf{C\textsubscript{4}}} & 
\multicolumn{1}{|p{2.71in}|}{\textbf{$ \ldots $ }} \\
\hhline{------}
%row no:2
\multicolumn{1}{|p{0.89in}}{P\textsubscript{1}} & 
\multicolumn{1}{|p{0.29in}}{\textbf{\textcolor[HTML]{FF0000}{0}}} & 
\multicolumn{1}{|p{0.29in}}{1} & 
\multicolumn{1}{|p{0.28in}}{0} & 
\multicolumn{1}{|p{0.29in}}{0} & 
\multicolumn{1}{|p{2.71in}|}{$ \ldots $ } \\
\hhline{------}
%row no:3
\multicolumn{1}{|p{0.89in}}{P\textsubscript{2}} & 
\multicolumn{1}{|p{0.29in}}{1} & 
\multicolumn{1}{|p{0.29in}}{\textbf{\textcolor[HTML]{FF0000}{1}}} & 
\multicolumn{1}{|p{0.28in}}{0} & 
\multicolumn{1}{|p{0.29in}}{1} & 
\multicolumn{1}{|p{2.71in}|}{$ \ldots $ } \\
\hhline{------}
%row no:4
\multicolumn{1}{|p{0.89in}}{P\textsubscript{3}} & 
\multicolumn{1}{|p{0.29in}}{0} & 
\multicolumn{1}{|p{0.29in}}{0} & 
\multicolumn{1}{|p{0.28in}}{\textbf{\textcolor[HTML]{FF0000}{1}}} & 
\multicolumn{1}{|p{0.29in}}{1} & 
\multicolumn{1}{|p{2.71in}|}{$ \ldots $ } \\
\hhline{------}
%row no:5
\multicolumn{1}{|p{0.89in}}{P\textsubscript{4}} & 
\multicolumn{1}{|p{0.29in}}{1} & 
\multicolumn{1}{|p{0.29in}}{0} & 
\multicolumn{1}{|p{0.28in}}{1} & 
\multicolumn{1}{|p{0.29in}}{\textbf{\textcolor[HTML]{FF0000}{0}}} & 
\multicolumn{1}{|p{2.71in}|}{$ \ldots $ } \\
\hhline{------}
%row no:6
\multicolumn{1}{|p{0.89in}}{P\textsubscript{5}} & 
\multicolumn{1}{|p{0.29in}}{0} & 
\multicolumn{1}{|p{0.29in}}{1} & 
\multicolumn{1}{|p{0.28in}}{0} & 
\multicolumn{1}{|p{0.29in}}{1} & 
\multicolumn{1}{|p{2.71in}|}{$ \ldots $ } \\
\hhline{------}
%row no:7
\multicolumn{1}{|p{0.89in}}{$ \ldots $ } & 
\multicolumn{1}{|p{0.29in}}{$ \ldots $ } & 
\multicolumn{1}{|p{0.29in}}{$ \ldots $ } & 
\multicolumn{1}{|p{0.28in}}{$ \ldots $ } & 
\multicolumn{1}{|p{0.29in}}{$ \ldots $ } & 
\multicolumn{1}{|p{2.71in}|}{$ \ldots $ } \\
\hhline{------}

\end{tabular}
 \end{table}

%%%%%%%%%%%%%%%%%%%% Table No: 1 ends here %%%%%%%%%%%%%%%%%%%%

\begin{justify}
Now, by using Cantor’s diagonalization logic, Alice is convinced that no matter what the outputs of P\textsubscript{1}, P\textsubscript{2}, P\textsubscript{3}, P\textsubscript{4} etc. are, Bob will never be able to cover all the choice of sequences that she has. That is because if she were to reverse the diagonal characters of the matrix above, that string can never have been outputted by Bob’s program. Hence, since Alice’s set of choices includes those that Bob’s computer program can never replicate, she concludes that it can never win the game if it were to be played for an arbitrarily long period of time. This means that at some point, Bob’s computer is going to guess wrongly. Let us assume that the length of the sequence then is  \( n \) . We will now prove that for this length, the Alice Bob example is actually NP complete. We will then analyse the methodology and insights derived from this exercise in the next section to understand the underlying basis for the result.
\end{justify}\par

\begin{justify}
For proving NP completeness, we will need to demonstrate the following three conditions:
\end{justify}\par

\begin{itemize}
	\item \textit{Condition A}: Validating the output can be done in polynomial time\par

	\item \textit{Condition B:} Arriving at a solution non-deterministically can be done in polynomial time\par

	\item \textit{Condition C}: A known NP complete program can be converted to the Alice Bob example
\end{itemize}\par

\begin{justify}
We will now proceed to resolve each of the above conditions.
\end{justify}\par

\begin{itemize}
	\item Resolution of \textit{Condition A}: Given the final output of Bob’s program, checking the result with Alice’s sequence is obviously a very simple process comprising of  \( n  \) steps (where  \( n \)  is the length of the string). We need to check whether each character of Bob’s sequence matches the corresponding character of Alice’s sequence. Hence, checking is an algorithm executable in polynomial time\par

	\item Resolution of \textit{Condition B}: Guessing a perfect solution non-deterministically is again just an  \( n \)  step process comprising of guessing each incremental character of the sequence \par

	\item Resolution of \textit{Condition C}: For this, we choose the  \( 3 \) - \( SAT \)  problem, which we know is NP complete.  We need to show that any input to the  \( 3 \) - \( SAT \)  problem can be converted as an input to the Alice Bob problem. We can do so by considering each disjunctive clause of the form  \( l_{1} \bigvee l_{2} \bigvee l_{3} \)  (where  \( l_{1} \) , \(  l_{2} \)  and  \( l_{3} \)  are the literals) in the  \( 3 \) - \( CNF \)  as a character corresponding to the Alice Bob problem with two possible values. Indeed, the Alice Bob example is just a disguised form of the  \( 3 \) - \( SAT \)  problem (as both of them have  \( 2^{n}  \) permutations). 
\end{itemize}\par

\begin{adjustwidth}{0.25in}{0.0in}
To understand why the  \( 3 \) - \( SAT \)  problem is in NP, let us contrast it with the  \( 2 \) - \( SAT \)  problem (which is in P). With  \( 2 \) - \( SAT \) , each individual clause of the form  \( l_{1} \bigvee l_{2} \)  can be expressed in the form of  \( a \Rightarrow b \)  This is because since every disjunctive clause in the  \( CNF \)  has to be true, it follows that  \(  \sim  l_{1} \Rightarrow  l_{2} \)  and  \(  \sim  \)   \( l_{2} \Rightarrow  l_{1} \) . Hence, this is completely deterministic because if  \( a \)  is true so is  \( b \) and if  \( b \)  is false so is  \( a \) . In other words, there is only one choice for each literal depending on the value of the other literal. So, to check if a given  \( 2 \) - \( CNF \)  formula is satisfiable, we only have to check if there are any conflicting conditions for each of the literals in all the clauses it appears. In case any such conflicts are detected for a particular literal, then the whole formula will be unsatisfiable. \par

\end{adjustwidth}

\begin{adjustwidth}{0.25in}{0.0in}
On the other hand, in  \( 3 \) - \( SAT \) , each clause takes the form of \(  a \Rightarrow  \)   \( b  \bigvee c \) .  Now, if we set  \( a \)  to be true then either  \( b \)  or  \( c \)  can be true, i.e. there are two possibilities. Therefore for a  \( 3 \) - \( SAT \)   \( CNF \)  formula of size  \( n \) , there will be  \( 2^{n} \)  possibilities. In general, for any  \( k \) - \( SAT \)  problem, there will be  \(  \left( k-1 \right) ^{n} \)  possibilities (where \(  n \)  is the number of disjunctive clauses in the  \( CNF \)  formula).\par

\end{adjustwidth}

\begin{justify}
Hence, we have proved Alice Bob problem to be NP complete. But, since for this sequence of length  \( n \) , Bob’s program failed to yield the correct result (in polynomial time), this implies that: 
\end{justify}\par

\begin{justify}
P!=NP 
\end{justify}\par

	\item \textbf{Analysis of the P-NP result}\par

\begin{justify}
With regard to the methodology adopted above, the one possible objection which could be raised is that we assumed that for Bob to definitely lose the game, it will have to be played for an arbitrarily long period of time. Yet, we set a definite limit to the length of the string  \( n \)  when he loses. We will address this concern here. The use of Cantor’s diagonalization technique for the proof alludes to a linkage with the cardinalities of natural and real number spaces. The key insight here is that the outputs generated by Bob’s algorithm have a one-to-one mapping to the natural number space while Alice is able to pick her sequences from the real number space. Herein, lies the fundamental chasm that all computer programs confront, which is unbridgeable. Indeed any problem where the result set can be in the real number space can never be solved in a deterministic manner by computers (which operate in natural number space) simply because the cardinality of the former is much larger than the cardinality of the latter. This again is not very hard to intuitively understand. The very essence of an irrational number (which is a part of the real number space) is that there is no pattern to the digits. As such, it is impossible to predict it in a deterministic manner. 
\end{justify}\par

\begin{justify}
Now, the question may be asked as to what about the cases where the output does not map to the real number space, but to the natural number space. To simulate this situation, let us replace Alice and Bob with Computer A and Computer B respectively. Now, both computers are constrained to operate in the natural number space. However, the question here is if this by itself fundamentally alters anything, the answer to which is negative. The constraint imposed on Computer A is merely a result of its physical limitations (such as memory and processing power), which can have arbitrary improvements with time. As such, the range of numbers that they can generate is not purely rational, but what we can term as $``$quasi-irrational$"$ \footnote{ For replicating a truly random sequences of numbers, we could assume the utilization of a device such as a Geiger counter }. What we mean by this is that an irrational number abruptly truncated to a certain length due to size constraints does not make it any more predictable. Hence, it will be the case that even though Computer A will generate such quasi-irrational numbers, since they don’t fit any pattern, it will be impossible for Computer B to guess them in a deterministic manner. Since a determinate algorithm can never be designed to predict a number that has inherently indeterminate characteristics, the P-NP result we derived above would hold for this case as well. 
\end{justify}\par

\begin{justify}
The other intuition we would like to draw out here is the connection with Gödel’s theorems. In fact, the precise logic we employed here can also be used to prove Gödel’s theorems as demonstrated by Gusfield [2014]. It is proved therein that there are valid functions which map to the real number space, which can never be replicated by a computer operating in the rational\footnote{ We will use the terms rational numbers space and natural numbers space interchangeably as their cardinalities are equal } number space. Hence, such functions are non-computable leading to the Gödel’s incompleteness result. This is directly analogous to the Alice Bob example here. In other words, the underlying basis of both problems is on account of the disparity in the cardinalities of the two spaces (i.e. real and rational). We have already demonstrated the causal connection of dual identities to the Gödel’s case. Hence, is it then possible that there is also a connection between computability, duality and the real and rational spaces? This is the question that we will tackle in the subsequent sections. By doing so, we will be able to explore all the nuances of P-NP problem fully (which will come together in the final section of this paper).
\end{justify}\par

	\item \textbf{Zeno’s paradox}\par

\begin{justify}
The dichotomy of rational and real numbers was best captured by Zeno, a Greek philosopher who lived between 490-430 BC, approximately 20 years before Socrates.  He devised a set of paradoxes, which together rule out the possibilities of space (and time) being either continuous or discrete. There are several paradoxes he constructed (Salmon [1970]), but we pick one from each category (i.e. continuous and discrete) for simplicity’s sake. 
\end{justify}\par

\begin{justify}
The so called $``$dichotomy paradox$"$  goes like this – suppose Homer wishes to walk to the end of a path. Before he can get there, he must get halfway there. Before he can get halfway there, he must get a quarter of the way there. Before traveling a quarter, he must travel one-eighth; before an eighth, one-sixteenth and so on.\ Since such divisions of intervals can be done infinite times, it will therefore require infinite tasks for Homer to reach his destination, which Zeno claims will be impossible for Homer to perform. Zeno uses this as an example to refute the continuity of space and time.  
\end{justify}\par

\begin{justify}
On the other hand, he uses the so called $``$arrow paradox$"$  to refute the discreteness of space and time. Herein, Zeno states that for motion to occur, a moving object such as an arrow in flight must change the position which it occupies. At any one instant of time, the arrow is where it is occupying a portion of space equal to itself. It is neither moving to where it is, nor to where it is not. It cannot move to where it is not, because no time elapses for it to move there; it cannot move to where it is, because it is already there. In other words, at every instant of time there is no motion occurring. If everything is motionless at every instant, and time is entirely composed of instants, then motion is impossible.
\end{justify}\par

\begin{justify}
Zeno’s paradoxes have elicited responses from mathematicians and philosophers alike. Since Zeno pre-dates the era of modern mathematics wherein concepts such as summation of infinite series and calculus were completely unknown, one could attempt to resolve this paradox strictly on mathematical terms by summing up the infinite series of: 
\end{justify}\par

\begin{justify}
 \( \frac{1}{2}+\frac{1}{4}+\frac{1}{8}+ \ldots  \) \tab \  
\end{justify}\par

\begin{justify}
which converges to 1. But, most philosophers have rejected such a strictly mathematical resolution of the paradox as it fails to satisfactorily address the basic issue of whether an infinite number of tasks, no matter how small can be accomplished in a finite amount of time. There have been some such as James [1948] and Whitehead [1948] who have questioned the very rationale behind representing time as a continuum as each observation itself requires a minimal length of time. Their arguments rely on the inadequacy of a continuous space for describing temporal processes. However, others such as Grünbaum [1968] have argued that strictly from a scientific standpoint, time is indeed infinitely divisible. Despite this, the need to build a mathematical construct to adequately address the paradoxes was undeniable resulting in an ongoing effort towards the same. Vlastos [1967] and McLaughlin and Miller [1992] have tried to build models based on mapping the full interval of space into finite number of sub-intervals based on whether an observable event has actually taken place in that particular interval. In this way, they circumvent the need to break up an interval into infinitely many intervals. Lee [1965] on the other hand, questions the very validity of Zeno’s methodology arguing that successively halving any finite interval of space only yields rational numbers and not irrational ones and therefore is unable to cover the entire continuity of space. 
\end{justify}\par

\begin{justify}
It is perhaps fair to say that none of the approaches hitherto adopted have satisfactorily addressed all the issues concerned. This is because neither have any of them made their way into mainstream mathematics nor have they engendered any active ongoing interest in the philosophical community to continue to probe them. Our approach herein would be to get to the heart of the matter from a theoretical standpoint and then build a mathematical construct to address the same, which we will do in the succeeding sections.
\end{justify}\par

	\item \textbf{Mapping of real numbers}\par

\begin{justify}
 The shortcomings of our current understanding of real numbers both from theoretical (Bell [2005], Ehrlich [2013]) as well as philosophical standpoints (Snapper [1979], Ferreirós [2008]) are well documented. Indeed, it has prompted some to call it the $``$most fundamental and important problem in mathematics$"$ \footnote{ Attributable to N J Wildberger who has a series of videos posted on the internet where he also calls the current\  theory of real numbers as a $``$logical sham$"$ . }. At the heart of this issue is – how does one actually define real numbers, and after having done that, what are the arithmetic principles governing them. One of the popular methodologies for doing so is based on the so called Dedekind cut (Dedekind [1963]), wherein the number line is partitioned into two parts – \textit{A} and \textit{B}, such that all elements of \textit{A} are less than all elements of \textit{B}, and \textit{A} contains no greatest element. An irrational number is defined by such a cut, wherein the number lies in neither of the two partitions. As an example of a Dedekind cut representing an irrational number, we may take the positive square root of 2. This can be defined by the set  A =  \(  \{ x \in Q:x<0 \vee x^{2}<2  \}  \) . The idea here is that any real number can be uniquely mapped using such sets of rational numbers engendered from such partitions. 
\end{justify}\par

\begin{justify}
Another such methodology for defining irrational numbers is to think of them as successive approximations of Cauchy sequences. E.g. – for  \(  \pi  \)  = 3.14159$ \ldots $  the Cauchy sequence would be:
\end{justify}\par

\begin{justify}
3, 3.1, 3.14, 3.141, 3.1415, 3.14159, $ \ldots $ .
\end{justify}\par

\begin{justify}
But, the question which has probably not been asked till now, which we do so here is – if such methods of defining an irrational number actually map it to a point on the number line. And the answer to that as we prove in the theorem below is - no.
\end{justify}\par

\begin{justify}
\textbf{\textit{Theorem.}} An irrational number cannot be represented as a point in the number line.
\end{justify}\par

\begin{justify}
\textbf{\textit{Proof. }}We set the initial stage by first representing all rational numbers as points in the number line. This is done by first mapping all the natural numbers, followed by integers and finally fractions based on well-known set theoretic methods. Since all of these numbers are well ordered, we can create a well-defined ordered mapping of them on the number line.  Using these mappings of rational numbers, we then attempt to map the rest of the real numbers, i.e. the irrational numbers. For this, we define an operation called the $``$point locator operation$"$   \(  \Gamma  \) , which helps map any real number  \( x  \) to a point on the number line.  \(  \Gamma  \)  will operate as follows. 
\end{justify}\par

\begin{justify}
We begin by choosing any arbitrary rational number  \(  \alpha _{1} \) , which maps to a definite point on the number line. If  \(  \alpha _{1}<x \) , we choose a point which corresponds to a number  \(  \alpha _{2} \ni  \alpha _{1}< \alpha _{2} \leq x \) . If  \(  \alpha _{2}=x \) ,  \(  \Gamma  \)  has been successful in mapping x to a point on the number line. Else, we repeat the process till the end goal is achieved.
\end{justify}\par

\begin{justify}
On the other hand, if  \(  \alpha _{1}>x \) , we choose a point which corresponds to a number  \(  \alpha _{2} \ni  \alpha _{1}> \alpha _{2} \geq x \) . If  \(  \alpha _{2}=x \) ,  \(  \Gamma  \)  has been successful in mapping x to a point on the number line. Else, we similarly repeat the process till the end goal is achieved. 
\end{justify}\par

\begin{justify}
If  \( x \)  is a rational number, it is trivial to prove that  \(  \Gamma  \)  will be immediately successful as it can just choose the point  \(  \alpha _{2} \)  that precisely maps to  \( x \)  itself.
\end{justify}\par

\begin{justify}
Next, consider the case where  \( x \)  is an irrational number and  \(  \alpha _{1}<x \) . Now,  \(  \Gamma  \)  will choose a point which corresponds to a number  \(  \alpha _{2} \ni   \alpha _{2}> \alpha _{1} \bigwedge  \alpha _{2} \leq x \) . But, since we haven’t mapped any irrational numbers yet, 
\end{justify}\par

 \(  \Rightarrow  \alpha _{2}  \) is also not irrational\par

 \(  \Rightarrow  \alpha _{2} \neq x  \) (as  \( x \)  is irrational)  \(  \Rightarrow  \alpha _{2}<x \) \par

\begin{justify}
Now, since there is always a rational number between any two numbers, we can choose a point which corresponds to a number  \(  \alpha _{3} \ni   \alpha _{3}> \alpha _{2} \bigwedge  \alpha _{3} \leq x \) 
\end{justify}\par

\begin{justify}
But for the same reason as before,  \(  \alpha _{3}<x \) . Hence, it follows that even if  \(  \Gamma  \)  operation is performed repeatedly, it can never map to a point that will exactly represent  \( x \)  on the number line. The same logic will hold true even if  \(  \alpha _{1}>x \) . 
\end{justify}\par

\begin{justify}
Hence, since  \(  \Gamma  \)  has failed to map even a single irrational number, it follows that none of the irrational numbers can be represented as points in the number line.
\end{justify}\par

\begin{justify}
\ \ \ \ \  The intuition behind the above theorem must be quite obvious. It is just an explicit way of stating that methods such as the Dedekind cut and Cauchy sequence are merely limiting and not precise definitions of irrational numbers. Nonetheless, this simple yet deep result has not been adequately absorbed, represented and formalized in the prevailing mathematical setup.\  The implication of this result is that we are left with contrasting requirements for rational and irrational numbers – the former being discrete and represented by points while the latter being continuous and represented by regions in space. We will proceed to build a mathematical framework for harmonising these requirements into a single space. But before that, we will briefly deal with the topic of Continuum Hypothesis in the next section.
\end{justify}\par

	\item \textbf{The Continuum Hypothesis}\par

\begin{justify}
A complete exposition on the Continuum Hypothesis (CH) is clearly out of scope for this paper and hence we will not attempt to do it. But, we will touch upon it here because it prepares a good ground for the mathematical framework that we will describe in the next section. CH is one of the most important unresolved problems in set theory. In 1874, Cantor had demonstrated the difference in cardinalities of the real and rational number spaces (Cantor [1874]). Cantor proved that while there is a one-to-one correspondence between the rational numbers and the natural numbers, there is no one-to-one correspondence between the real numbers and the natural numbers.  He thus concluded that the cardinality of the set of real numbers is larger than the cardinality of the set of rational numbers. But, he soon went on to seek an answer to the question of whether there were any other infinite sets of real numbers that were of intermediate size (Cantor [1878]). In other words, whether there was an infinite set of real numbers that could not be put into one-to-one correspondence with the natural numbers and could also not be put into one-to-one correspondence with the real numbers. The Continuum Hypothesis is simply the statement that there is no such set of real numbers. Not only was Cantor unable to resolve this question, but CH also remains an unsolved problem till date. In fact, it is probably unique in that there are proofs to demonstrate that it can neither be disproved (Gödel [1940]) nor proved (Cohen [1963], Cohen [1964]).
\end{justify}\par

\begin{justify}
For illustrating the possibility of such a space that has characteristics of both continuous and discrete spaces, consider the following sequence:
\end{justify}\par

\begin{justify}
0.111528111835111364111739111082111604111027$ \ldots $ .
\end{justify}\par

\begin{justify}
Do you notice anything here? On close inspection it can be seen that the sequence $``$111$"$  occurs regularly with random three digits interspersed between two such consecutive instances\ of 111’s.  If this pattern were to continue indefinitely, it is possible (indeed trivial) to predict the same with 50$\%$  accuracy. Actually, an even better way of articulating that would be to say that 50$\%$  of this space can be predicted with 100$\%$  accuracy (and hence computable), while 50$\%$  can be predicted with 0$\%$  accuracy (and hence non-computable). But the string as a whole is non-computable since it is irrational\footnote{ These could be termed as $``$pseudo-irrational$"$  numbers as they appear to be irrational, but on closer inspection, have some level of predictability. In contrast, $``$quasi-irrational$"$  numbers appear to be rational, but don’t have any predictability }. The question can then be asked if it is possible to treat such a sequence as a single combination of rational and real spaces and if so, what is the mathematics underlying it. The answer to that will be given in the next section.
\end{justify}\par

	\item \textbf{Defining the $``$Hybrid Space$"$ }\par

Let us begin by articulating the intuition behind the hybrid space, a proposed symbolic representation of space catering to the dual requirements of real numbers. Suppose  \( U \)  represents the number line of natural numbers and  \( V \) is a continuous space of positive real numbers, i.e.  \( V \subset  \)   \( R^{+} \) \par

\begin{justify}
Let us now define a mapping  \( f \)  from each discrete element in  \( U \)  to a continuous interval in  \( V \) 
\end{justify}\par

\begin{justify}
Hence,  \( f:u_{i} \rightarrow  \left[ v_{j},v_{k} \right]  where  \)   \( u_{i} \in U;~  \left[ v_{j},v_{k} \right]  \in V \) 
\end{justify}\par

\begin{justify}
Our goal with the mapping function is to faithfully map  \( U \)  to  \( V \) \textsubscript{. }However, since  \( V \)  is a continuous space, elements from  \( U \)  can only be mapped to small intervals and not directly to points in  \( V. \)  Hence, this mapping function  \( f \)  denotes the error in mapping elements from  \( U \)  to  \( V \) \textsubscript{.}Greater the interval of the mapping function, larger will be the error in mapping points onto  \( V \) 
\end{justify}\par

\begin{justify}
Let us now define  \( ++ \)  as the increment operator in  \( U \) , which will move any element in  \( U \)  to the next element there. This will correspondingly move the mapping in  \( V \)  to the next interval there. In other words, if  \( u_{i} \)  moves to  \( u_{t} \)  then [ \( v_{j},v_{k} \) ] will move to [ \( v_{l},v_{m} \) ].  
\end{justify}\par

\begin{justify}
Let us define  \( u_{t}-u_{i} \) as the $``$virtual distance$"$  and  \( v_{m}- v_{k} \)  as the $``$actual distance$"$  
\end{justify}\par

\begin{justify}
Hence, if we want to map the virtual distance in  \( U \) \textsubscript{ }to the actual distance in  \( V \) , the error will depend on the mapping function  \( f \) . We will use this intuition while defining the hybrid space  \( H \) 
\end{justify}\par

\begin{justify}
Now, a hybrid space  \( H \)  is one which has both discrete and continuous elements in the same space. Hence, unlike the above example where the virtual distance in  \( U \)  was mapped to the actual distance in  \( V \) ; in  \( H \) , both of them will be mapped to the same space. In this case, let us define the mapping function,  \(  \mu _{ij} \)  in the following manner:
\end{justify}\par

\begin{justify}
 \(  \mu _{ij}: \left( i,j \right)  \rightarrow x_{ij} \) , where  \( x_{ij},i,j \in H \) 
\end{justify}\par

\begin{justify}
Notice here that  \( x_{ij} \)  represents a continuous interval, while  \( i,j \)  represent discrete elements. The intuition here is that through the mapping function, we are trying to infuse a semantic context of distance to two disparate points (which lack this semantic context). As such,  \( x_{ij} \)  denotes the actual distance, while the ordered pair  \(  \left( i,j \right)  \)  denotes the virtual distance.
\end{justify}\par

\begin{justify}
We now need to define how the addition operator will work in this space. The addition operation should indicate what the actual distance covered will be for any two pairs of discrete points that are $``$added$"$  to each other. For this, assume two ordered pairs  \(  \left( i,j \right)  \)  and  \(  \left( k,l \right)  \)  which map to adjacent intervals  \( x_{ij} \)  and  \( x_{kl} \) . The addition operation is now defined as follows:
\end{justify}\par

\begin{justify}
  \( x_{ij} \)  +  \( x_{kl}=  \left( j-i \right) + \left( l-k \right)  \) 
\end{justify}\par

\begin{justify}
Now, if  \( k=j \)  it implies that the mapping from the virtual to actual distances is perfect without any errors (because  \( x_{ij} \)  and  \( x_{kl} \)  are meant to be adjacent intervals). 
\end{justify}\par

\begin{justify}
Hence, in such a case, 
\end{justify}\par

\begin{justify}
 \( x_{ij} \)  +  \( x_{kl}=l-i \) 
\end{justify}\par

\begin{justify}
i.e. the actual distance is just the difference of the two discrete numbers (i.e. the higher of the second pair minus the lower of the first pair)
\end{justify}\par

\begin{justify}
Else, if  \( k>j, \)  we will define a term called convexity $``$ \(  \complement " \) , which will capture the accuracy of mapping as follows:
\end{justify}\par

\begin{justify}
 \(  \complement  = \frac{ \left( l-i \right) }{ \left( k-j \right) } \) \tab 
\end{justify}\par

\begin{justify}
Note that  \(  \complement   \in  Q,  \)  and  \(   \complement  >1 \) 
\end{justify}\par

\begin{justify}
The virtual distance is given by  \(  \left( l-i \right)  \)  whereas the actual distance is given by  \(  \left( j-i \right) + \left( l-k \right)  \) 
\end{justify}\par

\begin{justify}
 \(  \therefore  \)  Actual Distance = Virtual Distance  \(  \left( 1-\frac{1}{C} \right)  \) 
\end{justify}\par

\begin{justify}
Hence, if convexity  \(  \complement  \)  is  \( \infty,  \Rightarrow  \)  we have been able to map the discrete space to continuous space with absolute accuracy. This also resolves the Zeno’s paradox kind of situations as follows. If the space concerned is completely determinate, then Homer does not need to subdivide the space into further subsections to cover the distance. He can travel right through it from start point to end point because the actual distance is \textit{defined} as just the difference between the end and start points. On the other hand, if  \(  \complement  \)   \(  \neq \infty, \)  then there will be an error that Homer will have to live with, which cannot be avoided no matter how he divides the distance.
\end{justify}\par

\begin{justify}
Let us now illustrate these workings for the example of the sequence considered in the earlier section, i.e.
\end{justify}\par

\begin{justify}
0.111528111835111364111739111082111604111027$ \ldots $ .
\end{justify}\par

\begin{justify}
Herein, two adjacent intervals will be of the form [111xxx], [111xxx]. When represented visually it can be depicted as follows:
\end{justify}\par

\par 
 
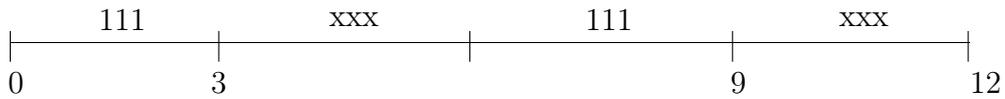
\begin{tikzpicture}

\path (0.19in,-0.48in) node [shape=rectangle,minimum height=0.26in,minimum width=0.3in,text width=0.27in,align=center]{0};

\path (1.24in,-0.48in) node [shape=rectangle,minimum height=0.26in,minimum width=0.27in,text width=0.24in,align=center]{3};

\path (5.21in,-0.48in) node [shape=rectangle,minimum height=0.26in,minimum width=0.53in,text width=0.48in,align=center]{12};

\path (3.93in,-0.48in) node [shape=rectangle,minimum height=0.26in,minimum width=0.28in,text width=0.25in,align=center]{9};

\draw (1.24in,-0.21in) -- (1.24in,-0.37in); 

\draw (0.16in,-0.27in) -- (5.13in,-0.27in); 

\path (4.58in,-0.15in) node [shape=rectangle,minimum height=0.26in,minimum width=0.53in,text width=0.48in,align=center]{xxx};

\path (3.26in,-0.15in) node [shape=rectangle,minimum height=0.26in,minimum width=0.53in,text width=0.48in,align=center]{111};

\path (1.94in,-0.15in) node [shape=rectangle,minimum height=0.26in,minimum width=0.53in,text width=0.48in,align=center]{xxx};

\path (0.74in,-0.15in) node [shape=rectangle,minimum height=0.26in,minimum width=0.53in,text width=0.48in,align=center]{111};

\draw (3.9in,-0.21in) -- (3.9in,-0.37in); 

\draw (2.54in,-0.21in) -- (2.54in,-0.37in); 

\draw (0.16in,-0.21in) -- (0.16in,-0.37in); 

\draw (5.12in,-0.21in) -- (5.12in,-0.37in); 
\setstretch{1.0}

\end{tikzpicture}
\begin{justify}
The convexity of the above space will be
\end{justify}\par

\setstretch{2.0}
\begin{justify}
 \( = \frac{ \left( 12-0 \right) }{ \left( 9-3 \right) }=2~  \) \tab 
\end{justify}\par

\begin{justify}
However, the advantage of representing the sequence in such a hybrid space is that a function whose range is in such a space is now computable. That is because the sequence of adjacent intervals of the form [111xxx] can be mapped to the natural number space whereas in its original form, the sequence was an irrational number and hence non-computable\footnote{ However, what we gain by computability, we compromise in accuracy }.\  
\end{justify}\par

\begin{justify}
We will also now highlight the significance of the convexity term  \(  \complement  \) . For this, consider another similar sequence to the one above but in this case, the intervals are of the form [111111xxx] (i.e. for every six 1’s there are three random digits). The convexity of such a space will be 3. Now, if we were to represent both these irrational numbers only in terms of regions of space, they could be done so as [0.11,0.12]. In other words, both of them can be represented by a space with lower bound 0.11 and upper bound 0.12. In doing so, one would have believed that they are both subject to the same level of uncertainty in mapping. However, as we have demonstrated, this fails to capture all the information we have about these two numbers (i.e. their convexities). Hence, a complete representation of an irrational number  \( x \)  will be of the form:
\end{justify}\par

\begin{justify}
  \( x\mapsto \left[ x_{i},x_{j}, \complement   \right] , \)  where  \( x_{i} \)  is the lower bound,  \( x_{j} \)  is the upper bound and  \(  \complement   \) is the convexity.
\end{justify}\par

\begin{justify}
Finally, we can also visualize a continuous space where different intervals have different convexities as follows: 
\end{justify}\par

\par 
 
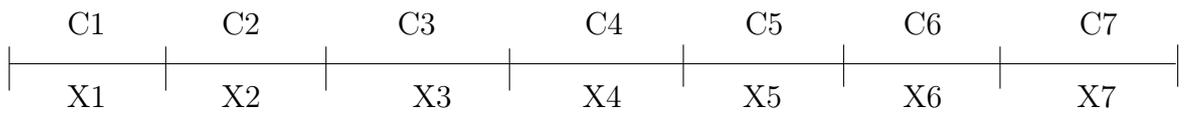
\begin{tikzpicture}

\draw (-0.03in,0.02in) -- (-0.03in,-0.21in); 

\path (5.6in,-0.24in) node [shape=rectangle,minimum height=0.27in,minimum width=0.36in,text width=0.32in,align=center]{X7};

\path (4.7in,-0.24in) node [shape=rectangle,minimum height=0.27in,minimum width=0.36in,text width=0.32in,align=center]{X6};

\path (3.87in,-0.24in) node [shape=rectangle,minimum height=0.27in,minimum width=0.36in,text width=0.32in,align=center]{X5};

\path (3.04in,-0.24in) node [shape=rectangle,minimum height=0.27in,minimum width=0.36in,text width=0.32in,align=center]{X4};

\path (2.16in,-0.24in) node [shape=rectangle,minimum height=0.27in,minimum width=0.36in,text width=0.32in,align=center]{X3};

\path (1.17in,-0.24in) node [shape=rectangle,minimum height=0.27in,minimum width=0.36in,text width=0.32in,align=center]{X2};

\path (0.37in,-0.24in) node [shape=rectangle,minimum height=0.27in,minimum width=0.36in,text width=0.32in,align=center]{X1};

\path (5.61in,0.14in) node [shape=rectangle,minimum height=0.27in,minimum width=0.36in,text width=0.32in,align=center]{C7};

\path (4.7in,0.14in) node [shape=rectangle,minimum height=0.27in,minimum width=0.36in,text width=0.32in,align=center]{C6};

\path (3.88in,0.14in) node [shape=rectangle,minimum height=0.27in,minimum width=0.36in,text width=0.32in,align=center]{C5};

\path (3.05in,0.14in) node [shape=rectangle,minimum height=0.27in,minimum width=0.36in,text width=0.32in,align=center]{C4};

\path (2.08in,0.14in) node [shape=rectangle,minimum height=0.27in,minimum width=0.36in,text width=0.32in,align=center]{C3};

\path (1.17in,0.14in) node [shape=rectangle,minimum height=0.27in,minimum width=0.36in,text width=0.32in,align=center]{C2};

\path (0.37in,0.14in) node [shape=rectangle,minimum height=0.27in,minimum width=0.36in,text width=0.32in,align=center]{C1};

\draw (5.1in,0.02in) -- (5.1in,-0.2in); 

\draw (6.02in,0.03in) -- (6.02in,-0.19in); 

\draw (4.29in,0.03in) -- (4.29in,-0.19in); 

\draw (3.46in,0.03in) -- (3.46in,-0.19in); 

\draw (2.56in,0.01in) -- (2.56in,-0.21in); 

\draw (1.61in,0.02in) -- (1.61in,-0.21in); 

\draw (0.78in,0.02in) -- (0.78in,-0.21in); 

\begin{scope}[yscale=1,xscale=-1,xshift=-5.98in]
	\draw (-0.03in,-0.07in) -- (6.01in,-0.07in); 
\end{scope}
\setstretch{1.0}

\end{tikzpicture}
\begin{justify}
(where C\textsubscript{1}, C\textsubscript{2}, C\textsubscript{3} etc. are the convexities of the virtual intervals x\textsubscript{1}, x\textsubscript{2}, x\textsubscript{3} etc.)
\end{justify}\par

\setstretch{2.0}
\begin{justify}
In such a case, the distance travelled will need to be computed as follows:
\end{justify}\par

\begin{justify}
If the actual distance of each interval is a constant, it follows that: 
\end{justify}\par

\begin{justify}
 \( x_{1} \left( 1-\frac{1}{C_{1}} \right) =x_{2} \left( 1-\frac{1}{C_{2}} \right) =x_{3} \left( 1-\frac{1}{C_{3}} \right)  \ldots   \)  
\end{justify}\par

From the above, we can deduce that the total actual distance travelled in such a scenario will be:\par

=  \( x_{1} \left[  \frac{ \left( 1-\frac{1}{C_{1}} \right) }{ \left( 1-\frac{1}{C_{2}} \right) } + \frac{ \left( 1-\frac{1}{C_{1}} \right) }{ \left( 1-\frac{1}{C_{3}} \right) } + \frac{ \left( 1-\frac{1}{C_{1}} \right) }{ \left( 1-\frac{1}{C_{4}} \right) }+ \ldots  \frac{ \left( 1-\frac{1}{C_{1}} \right) }{ \left( 1-\frac{1}{C_{n}} \right) }~  \right]  \) \par

=  \( x_{1} \left[ \frac{C_{1}-1}{C_{1}} \right]   \sum _{2}^{n}\frac{C_{i}}{C_{i}-1} \) \par

From our previous discussions in Section \(  2 \) , it will also be evident that for the 2-SAT problem, the convexity of solution space will be  \( \infty \) , while for the 3-SAT problem, it will be 3. This is so because it can be represented in the form [11x][11x]$ \ldots $   (as the introduction of the third literal in a clause causes the indeterminateness).\par

	\item \textbf{Philosophical implications of hybrid space}\par

\begin{justify}
We begin by stating the obvious, i.e.- that a number is meant to aid two kinds of activities – counting and measurement. These correspond to its existence in the discrete and continuous spaces respectively.  We know that in a discrete space, a number gets defined entirely by recursively executing a certain incremental rule as we do in a formal system like Peano arithmetic. So, we fix the origin 0 axiomatically and then define addition as a recursive function  \( F \) which maps two numbers to another one as follows:
\end{justify}\par

\begin{justify}
 \( a+0=a \) , \tab \tab \tab \tab (1)
\end{justify}\par

\begin{justify}
 \( a+F \left( b \right) =F \left( a+b \right)  \) \tab \tab (2)\tab 
\end{justify}\par

\begin{justify}
In\ this way, the entire discrete space gets defined in the form of incremental steps of the recursive operation.  Hence, there is no concept of a distance, but just a count of the number of operations that need to be performed to arrive at a number from the origin. Hitherto,\ such a discrete space has been treated as completely different from a continuous space, which is understandable since addition in the discrete space leads to an increase in plurality, i.e. 2 apples added to 3 apples give 5 distinct apples. In contrast, addition in continuous space leads to a decrease in plurality. Adding length 2 to length 3 gives a single length 5.  But, on the other hand, it is also the case that numbers cannot be defined in a mutually exclusively manner either in the discrete or continuous spaces. That is because for measuring length, we need a start point and end point in discrete space. Similarly, for defining a point in discrete space, we need to define how much distance is covered in each incremental operation. Hence, we need to define a single symbol which will satisfy both conditions. In other words, the number $``$5$"$  should stand for both 5 operations from the origin and also distance 5 from the origin. 
\end{justify}\par

\begin{justify}
But, as we demonstrated in the point mapping theorem, irrational numbers cannot be defined as discrete points and need to be represented by regions in space. As such, it is possible to view the concept of a number from the dual identity paradigm we have outlined earlier wherein, the states of a number can be defined as follows:
\end{justify}\par

\begin{itemize}
	\item \textit{State of Doing:} \textit{Separates} two points by a certain distance\par

	\item \textit{State of Being:} \textit{Is} a point which is a certain number of hops from the origin (as defined by a recursive increment function)
\end{itemize}\par

\begin{justify}
It merits mention here that when the underlying principles of quantum mechanics were first discovered, the notions of quantization and uncertainty associated with the physical world appeared completely counterintuitive. However, from our exposition here, such concerns should be completely laid to rest as these are not just constraints of the physical world, but a logical requirement of the very construct of space. 
\end{justify}\par

	\item \textbf{Conclusion}
\end{enumerate}\par

\begin{adjustwidth}{0.5in}{0.0in}
\begin{justify}
Let us tie in all the various points discussed so far.  A semantic understanding of Gödel’s theorems based on the duality of identity states of formal systems was developed.  However, that is just a part of the whole picture. To complete the understanding, this needs to be related to the dichotomy between real and rational number spaces. That will help us answer the question of why contradictions like Gödel’s incompleteness or Russell’s paradox can’t be avoided through some meta-level exclusion. It is not possible to do so for the following reasons:
\end{justify}\par

\end{adjustwidth}

\begin{adjustwidth}{0.75in}{0.0in}
a) All such cases involve definitions in first order logic, which are based on negation, e.g. - $``$A formal system which \textit{cannot} prove [certain propositions$ \ldots $ ]$"$  in the theorems of Gödel, or $"$ A barber who \textit{does} \textit{not} shave [a certain class of people$ \ldots $ ]$"$  in the paradox of Russell etc.\par

\end{adjustwidth}

\begin{adjustwidth}{0.75in}{0.0in}
b) Since the formal systems themselves are defined in rational number space, a negation of the same would map to the real number space (as the complement of rational number space is the irrational numbers space, which will be included in any negation)\par

\end{adjustwidth}

\begin{adjustwidth}{0.5in}{0.0in}
\begin{justify}
Hence, we cannot by any means exclude all possibilities of a negative definition through another meta-level definition (as the meta-level definition itself can only be constructed in the rational number space). In other words, since definitions map to the rational number space and negations map to the real number space, a negation can always be constructed, which cannot be annulled by another definition (either at first order or at higher orders). This is exactly the logic that is exploited in all diagonalization techniques a-la Cantor’s. Hence, such circularities based on negative definitions can never be broken.
\end{justify}\par

\end{adjustwidth}

\begin{adjustwidth}{0.5in}{0.0in}
\begin{justify}
This is where the demonstrated connection between the P-NP problem and the cardinalities of real and rational number spaces also assumes importance. The important point to note is that to solve the P-NP problem, we need to demonstrate that any polynomial time algorithm \textit{cannot} solve an NP complete problem. In other words, it is again a problem based on negation, alluding to a connection with the real number space. This explains the reliance on the diagonalization technique, which we employed to arrive at the proof. On the other hand, for reasons just explained earlier, an algorithm for identifying exclusions itself will be an NP problem as there will not be a deterministic way of solving it (apart from actually searching for solutions in the real number space). Hence, non-computability and NP are two facets of the same problem with a common connection to the real number space. To state in simple language - a problem becomes non-computable when the solution lies in the real number space, in which case the only way to solve it is to search for the solution through all possible alternatives.
\end{justify}\par

\end{adjustwidth}

\begin{adjustwidth}{0.5in}{0.0in}
\begin{justify}
But, as the point location theorem demonstrates, the representation of irrational numbers as points in space is a flawed one. As such, we conclude that the root cause of most unsolved problems in formal logic is a gap in the semantic understanding of the very concept of numbers. Specifically, this refers to the creation of two distinct definitions for rational and real numbers (i.e. discrete and continuous spaces) disjoint from each other. Hence, to address this, we have attempted to create a unified framework of numbers based on both rational and irrational numbers being mapped to a single space called the $``$hybrid space$"$ . In this notion of hybrid space, both discrete and continuous elements get \textit{defined} by a single parameter (called convexity) based on the degree of determinateness of that element. It is our belief that if the laws of the physical world are conceptualized on the basis of this space, many of the outstanding problems in theoretical physics also could get resolved. That will be the focus of our work going forward.
\end{justify}\par

\end{adjustwidth}

\vspace{\baselineskip}

 %%%%%%%%%%%%  Starting New Page here %%%%%%%%%%%%%%

\newpage

\vspace{\baselineskip}\begin{justify}
\textbf{REFERENCES}
\end{justify}\par

\begin{justify}
Bell, John Lane (2005), $``$The continuous and the infinitesimal in mathematics and philosophy$"$ , Polimetrica sas 
\end{justify}\par

\begin{justify}
Berto, Francesco (2011), There’s Something About Gödel: The Complete Guide to the Incompleteness Theorem, John Wiley $\&$  Sons
\end{justify}\par

\begin{justify}
Cantor, Georg (1874), $``$Über eine Eigenschaft des Inbegriffes aller reellen algebraischen Zahlen$"$ , Crelles Journal für reine und angewandte Mathematik 77: 258 – 262
\end{justify}\par

\begin{justify}
Cantor, Georg (1878), $``$Ein Beitrag zurMannigfaltigkeitslehre$"$ , Journal für reine und angewandte Mathematik\  84: 242–258
\end{justify}\par

\begin{justify}
Coffa, J. Alberto (1979), $``$The Humble Origins of Russell’s Paradox$"$ , Russell Nos. 33–34 
\end{justify}\par

\begin{justify}
Cohen, Paul (1963), "The Independence of the Continuum Hypothesis", Proc. Nat. Acad. Sci. U. S. A. 50: 1143-1148 
\end{justify}\par

\begin{justify}
Cohen, Paul (1964), "The Independence of the Continuum Hypothesis II", Proc. Nat. Acad. Sci. U. S. A. 51: 105-110 
\end{justify}\par

\begin{justify}
Cook, Stephen (1971), "The complexity of theorem proving procedures", Proceedings of the Third Annual ACM Symposium on Theory of Computing: 151–158 
\end{justify}\par

\begin{justify}
Dedekind, Richard (1963), $``$Continuity and Irrational Numbers$"$ , Essays on the Theory of Numbers, New York: Dover 
\end{justify}\par

\begin{justify}
Detlefsen, Michael (1979),\textit{ }$``$On interpreting Gödel's Second Theorem$"$ , Journal of Philosophical Logic 8: 297–313
\end{justify}\par

\begin{justify}
Ehrlich, Philip (2013), ed. $``$Real numbers, generalizations of the reals, and theories of continua Vol. 242$"$ , Springer Science $\&$  Business Media
\end{justify}\par

\begin{justify}
Ferreirós, José (2008), "The crisis in the foundations of mathematics", The Princeton companion to mathematics: 142-156
\end{justify}\par

\begin{justify}
Gödel, Kurt (1931), $``$On Formally Undecidable Propositions of Principia Mathematica and Related Systems$"$ , Translation of the German original Courier Corporation 
\end{justify}\par

\begin{justify}
Gödel, Kurt (1940), $``$The Consistency of the Continuum-Hypothesis$"$ , Princeton University Press
\end{justify}\par

\begin{justify}
Gusfield, Dan (2014), $``$Gödel for Goldilocks: A Rigorous, Streamlined Proof of (a variant of) Gödel’s First Incompleteness Theorem$"$ , arXiv:1409.5944
\end{justify}\par

\begin{justify}
Grünbaum, Adolf (1968), $``$Modern Science and Zeno’s Paradoxes$"$ , George Allen and Unwin Ltd, London 
\end{justify}\par

\begin{justify}
Hartmanis, Juris (1989), "Gödel, von Neumann, and the P = NP problem", Bulletin of the European Association for Theoretical Computer Science 38: 101–107
\end{justify}\par

\begin{justify}
James, William (1948), Some Problems of Philosophy, Longmans, Green $\&$  Co.,Ltd., London
\end{justify}\par

\begin{justify}
Kennedy, Juliette (ed.) (2014), Interpreting Gödel: Critical essays, Cambridge University Press 
\end{justify}\par

\begin{justify}
Lee, Harold (1965), "Are Zeno's Paradoxes Based on a Mistake?", Mind, New Series 74, no. 296: 563 – 570 
\end{justify}\par

\begin{justify}
McLaughlin, William and Miller, Sylvia (1992), $``$An Epistemological Use of Nonstandard Analysis to Answer Zeno’s Objections Against Motion$"$ , Synthese 92: 371-384
\end{justify}\par

\begin{justify}
Salmon, Wesley (ed.) (1970), $"$ Zeno’s Paradoxes$"$ , Bobbs-Merrill: 200—250
\end{justify}\par

\begin{justify}
Smith, Peter (2007)\textit{, }An Introduction to Gödel’s Theorems, Cambridge University Press 
\end{justify}\par

\begin{justify}
Snapper, Ernst (1979), "The three crises in mathematics: Logicism, intuitionism and formalism", Mathematics Magazine 52.4:\  207-216
\end{justify}\par

\begin{justify}
Vlastos, Gregory (1967), $``$Zeno of Elea$"$ , in Edwards, Paul (ed.), The Encyclopedia of Philosophy, Macmillan, New York, Vol. 8: 369–79
\end{justify}\par

\begin{justify}
Whitehead, Alfred North (1948), $``$Process and Reality$"$ , Macmillan, New York
\end{justify}\par

\vspace{\baselineskip}

\printbibliography
\end{document}